\documentclass[11pt]{amsart}

\usepackage{graphicx}
\usepackage{amscd}
\usepackage{amsmath}
\usepackage{amsfonts}
\usepackage{yfonts}
\usepackage{psfrag}
\usepackage{amssymb}
\usepackage{verbatim}
\usepackage{comment}
\newtheorem{theorem}{Theorem}[section]
\theoremstyle{plain}

\newtheorem{corollary}[theorem]{Corollary}

\newtheorem{lemma}[theorem]{Lemma}

\newtheorem{prop}[theorem]{Proposition}

\def\diam{{\rm diam}}
\def\dist{{\rm dist}}
\def\vphi{\varphi}

\def\Hk{{\mathcal H}}
\def\Ok{{\mathcal O}}

\def\sms{\setminus}
\def\th{\theta}
\def\Th{\Theta}

\newcommand{\lam}{\lambda}
\def\Lam{\Lambda}

\newcommand{\om}{\omega}
\def\Om{\Omega}
\newcommand{\Sig}{\Sigma}

\newcommand{\Gam}{\Gamma}
\newcommand{\sig}{\sigma}

\newcommand{\bi}{{\bf i}}

\newcommand{\R}{{\mathbb R}}
\newcommand{\Q}{{\mathbb Q}}

\newcommand{\Nat}{{\mathbb N}}
\newcommand{\Leb}{{\Hk^1}}
\def\Ak{{\mathcal A}}
\def\Sk{{\mathcal S}}
\def\Lk{{\mathcal L}}
\def\Kk{{\mathcal K}}
\def\Pk{{\mathcal P}}

\newcommand{\Wk}{{\mathcal W}}
\newcommand{\eps}{{\varepsilon}}
\newcommand{\es}{\emptyset}

\def\ov{\overline}

\begin{document}

\title[Visibility for self-similar sets]{Visibility for self-similar sets
of dimension one in the plane}

\author{K\'{a}roly Simon}
\address{K\'{a}roly Simon, Institute of Mathematics, Technical
University of Budapest, H-1529 B.O.box 91, Hungary}
\email{simonk@math.bme.hu}

\author{Boris Solomyak}
\address{Boris Solomyak, Box 354350, Department of Mathematics,
University of Washington, Seattle WA 98195}
\email{solomyak@math.washington.edu}

\date{\today}

 \thanks{
The research of K. S. was supported in part by OTKA Foundation grant T42496.
The research of B. S.  was supported in part by NSF grant DMS-0355187.
This collaboration was supported by NSF-MTA-OTKA
grant \#77.
}

\begin{abstract} We prove that a purely unrectifiable self-similar set
of finite 1-dimensional Hausdorff measure in the plane, satisfying the
Open Set Condition, has radial projection of zero length from
every point.
\end{abstract}

\maketitle

\thispagestyle{empty}

\section{Introduction}

For $a\in \R^2$, let $P_a$ be the radial projection from $a$:
$$
P_a:\ \R^2 \sms \{a\} \longrightarrow S^1,\ \ \ P_a(x) =
\frac{(x-a)}{|x-a|}\,.
$$
\begin{figure}[h!]
\psfrag{L}{$\Lambda  $}%
\psfrag{PaL}{$P_a(\Lambda ) $}%
\psfrag{a}{$a$}
\includegraphics*[bb=13 14 356 161,keepaspectratio,width=12cm]
{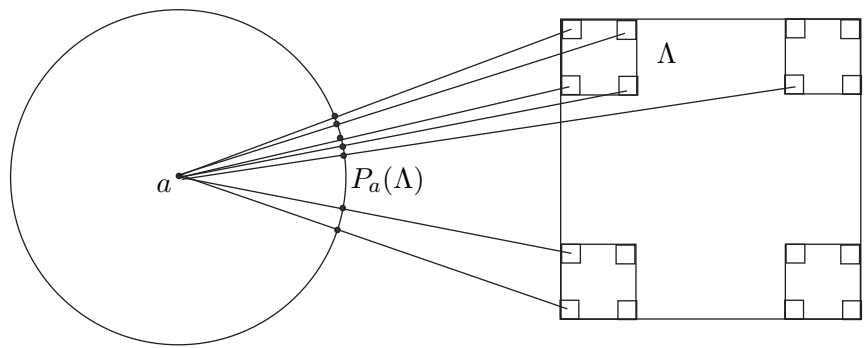} \caption{The radial projection of the four corner set}
\label{zeroth}
\end{figure}

A special case of our theorem asserts that the ``four corner Cantor
set" of contraction ratio $1/4$ has radial projection of zero length
from all points $a\in \mathbb{R}^2$. See Figure \ref{zeroth} where
we show the second-level approximation of  the four corner Cantor set
and the radial projection of some of its points.

Denote by $\Hk^1$ the one-dimensional Hausdorff measure. 
A Borel set $\Lam$ is a 1-{\em set} if $0 < \Hk^1(\Lam) <
\infty$. It is said to be {\em invisible from $a$} if $P_a(\Lam\sms\{a\})$
has zero length.

\begin{theorem} \label{th-main}
Let $\Lam$ be a self-similar 1-set in $\R^2$ satisfying
the Open Set Condition, which is not on a line.
Then $\Lam$ is invisible from every $a\in \R^2$.
\end{theorem}
                                                                                
Recall that a nonempty compact $\Lam$ is self-similar if
$\Lam = \bigcup_{i=1}^m S_i(\Lam)$ for some contracting similitudes $S_i$.
This means that
$$
S_i(x) = \lam_i\Ok_i x + b_i,
$$
where $0<\lambda_i<1$, $\Ok_i$ is an orthogonal transformation
of the plane, and $b_i
\in \R^2$. The Open Set Condition holds if there
exists an open set $V\ne \es$ such that $S_i(V) \subset V$ for all $i$ and
$S_i(V)\cap S_j(V)=\es$ for all $i\ne j$.
For a self-similar set satisfying the Open Set Condition, being a 1-set
is equivalent to $\sum_{i=1}^m \lam_i =1$.

A Borel set $\Lam$ is {\em purely unrectifiable} (or {\em 
irregular}), if $\Hk^1(\Lam\cap \Gam) = 0$ for every rectifiable curve $\Gam$.
A set $\Lam$ satisfying the assumptions of Theorem~\ref{th-main} is
purely unrectifiable by Hutchinson \cite{Hutch} (see also \cite{Mat2}).
A classical theorem of
Besicovitch \cite{besi} (see also \cite[Theorem 6.13]{falc}) says that a purely
unrectifiable 1-set has orthogonal projections of zero length on almost
every line through the origin. We use it in our proof.
                                                                                
In \cite[Problem 12]{matsurv} (see also \cite[10.12]{mattila})
Mattila raised the following question: Let $\Lam$ be a Borel set in
$\R^2$ with $\Hk^1(\Lam) < \infty$. Is it true that for $\Hk^1$ almost all
$a\in A$, the intersection $\Lam \cap L$ is a finite set for almost all
lines $L$ through $a$? If $\Lam$ is purely unrectifiable, is it true
that $\Lam \cap L = \{a\}$ for almost all lines through $a$?
Our theorem implies a positive answer
for a purely unrectifiable self-similar 1-set $\Lam$
satisfying the Open Set Condition.
The general case of a purely unrectifiable set remains open.
On the other hand, M. Cs\"ornyei and D. Preiss proved recently
that the answer to the first part of the question is negative
[personal communication].

Note that we prove a stronger property for our class of sets, namely, that
the set is invisible from {\em every} point $a\in \R^2$. It is easy
to construct examples of non-self-similar purely unrectifiable
1-sets for which this 
property fails. Marstrand \cite{mars} has an example of a
purely unrectifiable 1-set which is visible from a set of dimension one.

We do not discuss here other results and problems related to
visibility; see \cite[Section 6]{mattila} for a recent survey.
We only mention a result of Mattila \cite[Th.5.1]{mattila2}: 
if a set $\Lam$ has projections of zero length
on almost every line (which could have $\Hk^1(\Lam) = \infty$), then the
set of points $\Xi$ from which $\Lam$ is visible is a purely unrectifiable set
of zero 1-capacity. A different proof of this and a characterization of
such sets $\Xi$ is due to Cs\"ornyei \cite{marianna}.


\section{Preliminaries}

We have $S_i(x):=\lambda _i\mathcal{O}_ix+b_i$, where
$0<\lambda_i<1$,
$$\mathcal{O}_i=\left[\begin{array}{cr}
 \cos(\varphi_i) & -\varepsilon_i\sin(\varphi_i) \\
 \sin(\varphi_i) & \varepsilon_i\cos(\varphi_i)  \\
\end{array}
\right],$$  $\varphi _i\in [0,2\pi )$, and $\varepsilon _i\in
\left\{-1,1\right\}$ shows whether $\mathcal{O}_i$ is a rotation through the
angle $\varphi _i$ or a reflection about the line through the origin
making the angle $\varphi_i /2$ with the $x$-axis.

Let $\Sigma :=\left\{1,\dots ,m\right\}^\mathbb{N}$ be the symbolic
space. The natural projection $\Pi:\,\Sig \to \Lam$ is defined by
\begin{equation} \label{eq-nat}
\Pi(\bi) = \lim_{n\to \infty} S_{i_1\ldots i_n} (x_0), \ \ \ \mbox{where}\ \
\bi = (i_1 i_2 i_3\ldots) \in \Sig,
\end{equation}
and $S_{i_1\ldots i_n} = S_{i_1} \circ \cdots \circ S_{i_n}$.
The limit in (\ref{eq-nat}) exists and does not depend on $x_0$.
Denote $\lam_{i_1\dots i_n} = \lam_{i_1}\cdots \lam_{i_n}$ and
$\eps_{i_1\ldots i_k} = \eps_{i_1} \cdots \eps_{i_k}$.
We can write
$$
S_{i_1\dots i_n}(x)=\lambda _{i_1\dots i_n}\Ok_{i_1\dots
i_n}x+ b_{i_1\dots i_n},
$$
where
$$
\Ok_{i_1\dots i_n}:= \Ok_{i_1} \circ \cdots \circ \Ok_{i_n} =
\left[\begin{array}{cr}
 \cos(\varphi _{i_1\dots i_n}) & -\eps_{i_1\ldots i_n}\sin(\varphi _{i_1\dots i_n}) \\
 \sin(\varphi _{i_1\dots i_n}) & \eps_{i_1\ldots i_n}\cos(\varphi _{i_1\dots i_n})  \\
\end{array}\right],
$$
$$
\vphi_{i_1\ldots i_n} := \vphi_{i_1} + \eps_{i_1} \vphi_{i_2} +
\eps_{i_1 i_2} \vphi_{i_3} + \cdots + \eps_{i_1\ldots i_{n-1}}\vphi_{i_n},
$$
and
$$
b_{i_1\dots i_n} = b_{i_1} + \lam_{i_1} \Ok_{i_1} b_{i_2} + \cdots +
\lam_{i_1\ldots i_{n-1}} \Ok_{i_1\ldots i_{n-1}} b_{i_n}.
$$
\begin{sloppypar}
\noindent
Since $\sum_{i=1}^m \lam_i =1$, we can consider the probability product measure
$\mu = (\lam_1,\ldots,\lam_m)^{\Nat}$ on the symbolic space $\Sig$
and define the {\em natural measure} on $\Lam$:
$$
\nu = \mu \circ \Pi^{-1}.
$$
By a result of Hutchinson \cite[Theorem 5.3.1(iii)]{Hutch}, as a consequence
of the Open Set Condition we have
\begin{equation} \label{name}
\nu = c\Hk^1|_\Lam,\ \ \ \mbox{where}\ \ c = (\Hk^1(\Lam))^{-1}.
\end{equation}
\end{sloppypar}

To $\th\in [0,\pi)$ we associate the unit vector $e_\th = (\cos\th,\sin\th)$,
the line $L_\th = \{te_\th:\ t\in \R\}$, and the orthogonal projection
onto $L_\th$ given by $x\mapsto (e_\th\cdot x)e_\th$. It is more convenient
to work with the signed distance of the projection to the origin, which
we denote by $p_\th$:
$$
p_\th:\, \R^2\to \R,\ \ \ p_\th x = e_\th\cdot x.
$$
Denote $\Ak:= \{1,\ldots,m\}$ and let $\Ak^* = \bigcup_{i=1}^\infty \Ak^i$
be the set of all finite words over the alphabet $\Ak$.
For $u = u_1\ldots u_k \in \Ak^k$ we define the corresponding ``symbolic''
cylinder set by
$$
[u] = [u_1\ldots u_k] := \{\bi \in \Sig:\ i_\ell = u_\ell, \, 1 \le \ell \le k
\}.
$$
We also let 
$$
\Lam_u = S_u(\Lam) = \lam_u \Ok_u \Lam + b_u
$$
and call $\Lam_u$ the cylinder set of $\Lam$ corresponding to the word $u$.
Let $d_\Lam$ be the diameter of $\Lam$; then $\diam(\Lam_u) = \lam_u d_\Lam$.
For $\rho>0$ consider the ``cut-set''
$$
\Wk(\rho) = \{u\in \Ak^*:\ \lam_u\le \rho,\ \lam_{u'} >\rho\}
$$
where $u'$
is obtained from $u$ by deleting the last symbol. 
Observe that for every $\rho>0$,
\begin{equation} \label{decom}
\Lam = \bigcup_{u\in \Wk(\rho)} \Lam_u.
\end{equation}
In view of (\ref{name}), we have $\nu(\Lam_u \cap \Lam_v) = 0$ for distinct
$u,v\in \Wk(\rho)$, hence
$$
\nu(\Lam_u) = \lam_u\ \ \ \mbox{for all}\ u\in \Ak^*.
$$
Denote $\lam_{\min} := \min\{\lam_i:\ i\le m\}$; then
$\mu(\Lam_u) = \lam_u \in (\rho \lam_{\min},\rho]$ for $u\in \Wk(\rho)$.

We identify the unit circle $S^1$ with $[0,2\pi)$ and use additive notation
$\theta_1 + \theta_2$ understood mod $2\pi$ for points on the circle.
For  a Radon measure $\eta$ on the line or on $S^1$,
the upper density of $\eta$ with respect to $\Hk^1$ is defined by
$$
\ov{D}(\eta,t) = \limsup_{r\to 0} \frac{\eta([t-r,t+r])}{2r}\,.
$$
The open ball of radius $r$ centered at $x$ is denoted by $B(x,r)$.


\section{Proof of the main theorem}

In the proof of Theorem~\ref{th-main}
we can assume, without loss of generality, that $a\not\in
\Lam$, and
\begin{equation} \label{ass2}
\mbox{$P_a(\Lam)$ is contained in an arc of length less than $\pi$.}
\end{equation}
Indeed, $\Lam\sms \{a\}$ can be written as
a countable union of self-similar sets
$\Lam_u$ for $u \in \Ak^*$, of arbitrarily small diameter. If each of them
is invisible from $a$, then $\Lam$ is invisible from $a$.

Let
$$
\Om:= \{\bi \in \Sig:\ \forall\,u\in \Ak^*\ \exists\,n\ \mbox{such that}\
 \sig^n\bi \in [u]\},
$$
that is, $\Om$ is the set of sequences which contain each finite word over
the alphabet $\Ak= \{1,\ldots,m\}$. It is clear that every $\bi \in \Om$
contains each finite word infinitely many times and $\mu(\Sig\sms \Om) =0$.

\begin{lemma}[Recurrence Lemma] \label{lem-rec}
For every $\bi \in \Om,\ \delta>0$, and $j_1,\ldots,j_k \in \{1,\ldots,m\}$,
there are infinitely many $n\in \Nat$ such that
\begin{equation} \label{eq-rec}
\phi_{i_1\ldots i_n} \in [0,\delta],\ \eps_{i_1\ldots,i_n} =1,\ \ \mbox{and}
\ \ \sig^n\bi \in [j_1\ldots j_k].
\end{equation}
\end{lemma}

If the similitudes have no rotations or reflections, that is, $\phi_i = 0$
and $\eps_i = 1$
for all $i\le m$ (as in the case of the four corner Cantor set), 
then the conditions on $\phi$ and $\eps$ in (\ref{eq-rec}) hold
automatically and
the lemma is true by the definition of $\Omega$. The proof in the
general case is not difficult, but requires a detailed
case analysis, so we postpone
it to the next  section.

Let
$$
\Th := \{\theta \in [0,\pi):\ \Leb(p_\th(\Lam))=0\}\ \ \ \mbox{and}\ \ \
\Th' := (\Th + \pi/2) \cup (\Th + 3\pi/2)
$$
(recall that addition is considered mod $2\pi$).
Since $\Lam$ is purely unrectifiable, $\Leb([0,\pi)\sms \Th') = 0$
by Besicovitch's Theorem \cite{besi}.
The following proposition is the key step of the proof.
We need the following measures:
$$
\nu_a := \nu\circ P_a^{-1}\ \ \ \mbox{and}\ \ \ \nu_\th:= \nu\circ p_\th^{-1},
\ \theta \in [0,\pi).
$$
We also denote $ \Lam' = \Pi(\Om)$.

\begin{prop} \label{prop-dens}
If $\theta' \in P_a(\Lam') \cap \Th'$, then
$\ov{D}(\nu_a,\theta') = \infty$.
\end{prop}

\begin{sloppypar}
{\em Proof of Theorem~\ref{th-main} assuming Proposition~\ref{prop-dens}.}
By Proposition~\ref{prop-dens} and \cite[Lemma 2.13]{mattila} (a
corollary of the Vitali covering theorem), we obtain that
$\Leb(P_a(\Lam') \cap \Th') = 0$. As noted above,
$\Th'$ has full $\Leb$ measure in $S^1$. On the other hand,
$$
\mu(\Sig\sms \Om) = 0\ \Rightarrow\ \nu(\Lam \sms \Lam') = 0
\ \Rightarrow\ \Hk^1(\Lam \sms \Lam') = 0
\ \Rightarrow\ \Leb(P_a(\Lam \sms \Lam')) = 0,
$$
and we conclude that $\Leb(P_a(\Lam))=0$, as desired. \qed
\end{sloppypar} 

\medskip

{\em Proof of Proposition~\ref{prop-dens}.}
Let $x\in \Lam'$ and $\theta' = P_a(x) \in \Th'$. Let
$\th := \th' - \pi/2$ mod $[0,\pi)$. By the definition of $\Th'$ we have
$\Leb(p_{\th}(\Lam)) = 0$.

\smallskip

First we sketch the idea of the proof. Since $\Leb(p_{\th}(\Lam)) = 0$, we
have $\nu_\th\perp \Leb$, and this implies that for every $N\in \Nat$ there
exist $N$ cylinders of $\Lam$ approximately the same diameter (say,
$\sim r$), such that their projections to $L_\th$ are $r$-close to each other.
Then there is a line parallel to the segment $[a,x]$, whose
$Cr$-neighborhood contains all $\Lam_{u_j}, j=1,\ldots,N$.
By the definition of $\Lam' = \Pi(\Om)$, we can find similar copies of this
picture near $x\in \Lam'$ at arbitrarily small scales. The Recurrence Lemma
\ref{lem-rec} guarantees that these copies can be chosen with a small
relative rotation. This will give $N$ cylinders of $\Lam$
of diameter $\sim r_0 r$ contained in a $C'r_0r$-neighborhood of the
ray obtained by extending $[a,x]$. Since $a$ is assumed to be separated
from $\Lam$, we will conclude that
$\ov{D}(\nu_a,\th') \ge C''N$, and the proposition will follow.
Now we make this precise. The proof is illustrated in Figure 2.

\begin{figure}[h!]
\psfrag{t}{$\theta $}%
\psfrag{L}{$L_\theta $}%
\psfrag{lue}{$\Lambda _{u^{(1)}}$}%
\psfrag{luN}{$\Lambda _{u^{(N)}}$}%
\psfrag{lao}{$\Lambda _w $}%
\psfrag{lve}{$\Lambda _{v^{(1)}}$}%
\psfrag{lvN}{$\Lambda _{v^{(N)}}$}%
\psfrag{x}{$x$}%
\psfrag{z}{zoom in}%
\psfrag{a}{$a$}%
\psfrag{or}{$(0,0)$}%
\psfrag{pax}{$P_a(x)$}%
\psfrag{luk}{$\Lambda _{u^{(2)}}$}%
\psfrag{lvk}{$\Lambda _{v^{(2)}}$}%
\psfrag{d}{$\vdots$}%
\includegraphics*[bb=3 63 618 619,keepaspectratio,width=15cm]
{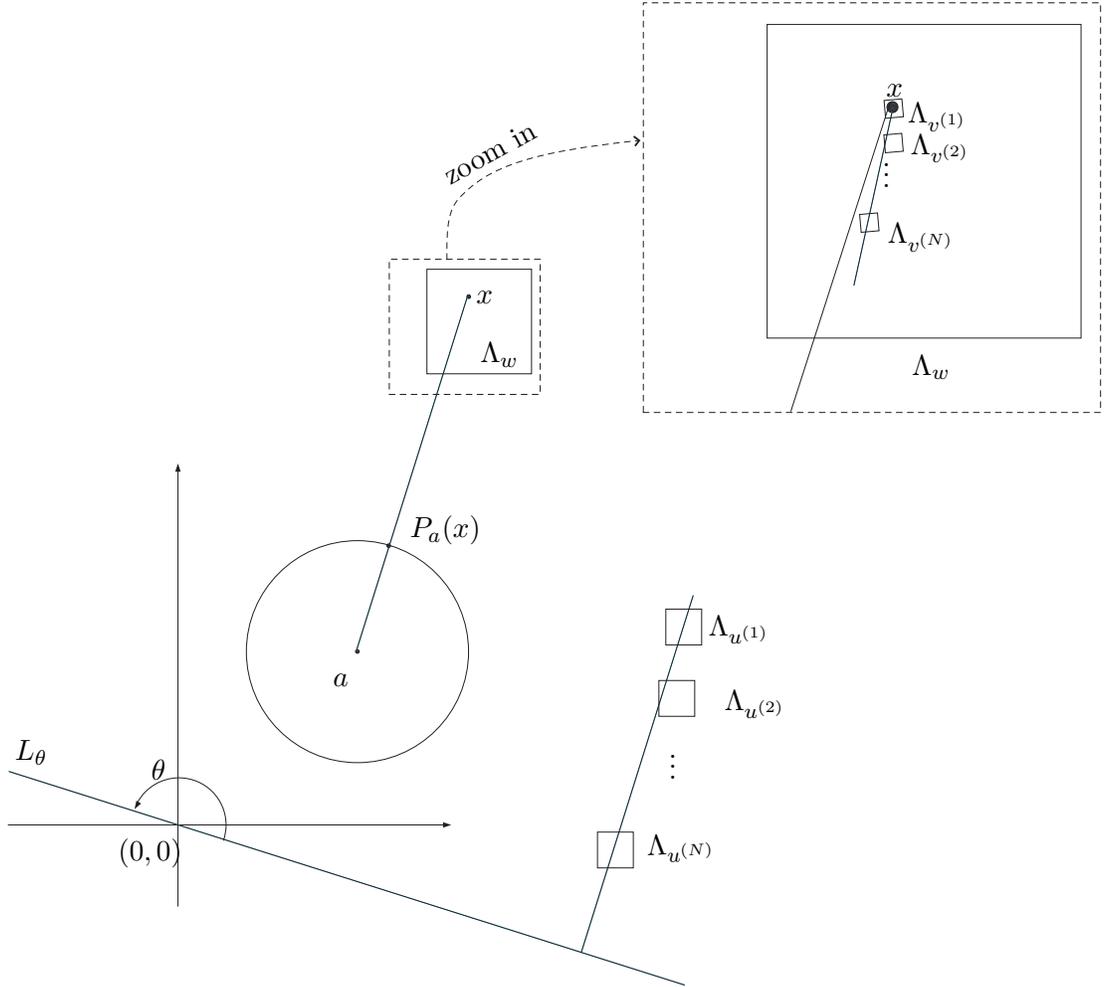} \caption{The cylinders of $\Lambda $ causing high
density} \label{first}
\end{figure}

\smallskip

{\sc Claim.}
{\em For each $N\in \Nat$ there exists $r>0$
and distinct $u^{(1)},\ldots, u^{(N)} \in \Wk(r)$ such that}
\begin{equation} \label{eq-pro1}
|p_{\th}(b_{u^{(j)}}-b_{u^{(i)}})| \le r,\ \ \forall\,i,j\le N.
\end{equation}
Indeed, for every $u\in \Ak^*$,
$$
\Lam_u = \lam_u \Ok_u \Lam + b_u\ \Rightarrow\
\Lam_u \subset B(b_u,d_\Lam \lam_u),
$$
hence for every interval $I\subset \R$ and $r>0$,
$$
\nu_{\th} (I) \le \sum_{u\in \Wk(r)}\{\lam_u:\ \dist(p_\th(b_u),I)
\le d_\Lam r\}.
$$
If the claim does not hold, then there exists $N\in \Nat$ such that
for every $t\in \R$ and $r>0$,
$$
\nu_{\th} ([t-r,t+r]) \le N(2(1+d_\Lam)+1)r.
$$
Then $\nu_{\th}$ is absolutely continuous with respect to $\Leb$,
which is a contradiction. The claim is verified. \qed

\smallskip

We are given that $x \in \Lam' = \Pi(\Om)$, which means that
$x = \pi(\bi)$ for an infinite sequence $\bi$ containing all finite words.
We fix $N\in \Nat$ and find $r>0$, $u^{(1)},\ldots, u^{(N)} \in \Wk(r)$ from
the Claim.
Then we apply Recurrence Lemma~\ref{lem-rec} with $j_1\ldots j_k:= u^{(1)}$
and $\delta = r$ to obtain infinitely many $n\in \Nat$ satisfying
(\ref{eq-rec}). Fix such an $n$. Denote
$$
w:= i_1\ldots i_n\ \ \ \mbox{and}\ \ v^{(j)}=wu^{(j)},\ j=1,\ldots,N.
$$
Observe that $\bi$ starts with $v^{(1)}$, so $x=\Pi(\bi) \in \Lam_{v^{(1)}}$,
hence
\begin{equation} \label{eq-pro2}
|p_{\th}(x-b_{v^{(1)}})|\le |x-b_{v^{(1)}}|
\le d_\Lam \lam_{v^{(1)}} \le d_\Lam \lam_w r.
\end{equation}
Here we used that $u^{(1)} \in \Wk(r)$, so $\lam_{v^{(1)}} =
\lam_w\lam_{u^{(1)}}\le \lam_w r$.
We have for $z\in \R^2$,
$$
\lam_{v^{(j)}} \Ok_{v^{(j)}} z + b_{v^{(j)}} = S_{v^{(j)}}(z) =
S_w\circ S_{u^{(j)}}(z) =
\lam_w\Ok_w(\lam_{u^{(j)}} \Ok_{u^{(j)}} z + b_{u^{(j)}}) + b_w,
$$
hence
$$
b_{v^{(j)}} = \lam_w\Ok_w b_{u^{(j)}} + b_w.
$$
It follows that
$$
p_{\th}(b_{v^{(i)}} - b_{v^{(j)}}) = \lam_w p_{\th} \Ok_w (b_{u^{(i)}}-
b_{u^{(j)}}).
$$
By (\ref{eq-rec}), we have $\eps_w=1$ and $\phi:= \phi_w \in [0,r)$;
therefore, $\Ok_w=R_\th$ is the rotation through the angle $\phi$.
One can check that $p_\th R_\phi = p_{\th-\phi}$,
which yields
\begin{equation} \label{eq-pro3}
|p_{\th}(b_{v^{(i)}} - b_{v^{(j)}})|= \lam_w |p_{\th-\phi}(b_{u^{(i)}}-
b_{u^{(j)}})|.
\end{equation}
Clearly, $\|p_\th - p_{\th-\phi}\|\le |\phi|\le r$, where $\|\cdot\|$
is the operator norm,
so we obtain from (\ref{eq-pro1}) and (\ref{eq-pro3}) that
$$
|p_{\th}(b_{v^{(i)}} - b_{v^{(j)}})|\le \lam_w(|b_{u^{(i)}}-b_{u^{(j)}}|r + r)
\le\lam_w(d_\Lam+1)r.
$$
Recall that $\bi$ starts with $v_1$, so $x=\Pi(\bi) \in \Lam_{v^{(1)}}$, hence
for each $j\le N$, for every $y\in \Lam_{v^{(j)}}$,
\begin{eqnarray}
|p_{\th}(x-y)| & \le & |x-b_{v^{(1)}}|+|p_{\th}(b_{v^{(1)}}- b_{v^{(j)}})| +
|b_{v^{(j)}}-y| \nonumber \\
& \le & d_\Lam (\lam_{v^{(1)}} +\lam_{v^{(j)}}) + \lam_w (d_\Lam+1)r\le
\lam_w (3d_\Lam + 1) r. \label{eq-pro4}
\end{eqnarray}
Now we need a simple geometric fact: given that
$$
P_a(x)= \th',\ \ \th = \th' + \pi/2\ {\rm mod}\ [0,\pi),\ \
|p_\th(x-y)| \le \rho,\ \ |y-a| \ge c_1,\ \ \mbox{and (\ref{ass2}) holds},
$$
we have
$$
|P_a(y)-\th'| = |P_a(y) -P_a(x)| =
\arcsin\frac{|p_\th(y-x)|}{|y-a|} \le
\frac{\pi}{2c_1} \rho.
$$
This implies, in view of (\ref{eq-pro4}), that for $c_2=
\pi(3d_\Lam + 1)/(2c_1)$,
$$
\nu_a([\th'-c_2 \lam_w r, \th'+c_2 \lam_w r]) \ge \sum_{j=1}^N
\nu(\Lam_{v^{(j)}})
= \sum_{j=1}^N \lam_{v^{(j)}} = \lam_w \sum_{j=1}^N \lam_{u^{(j)}} \ge
 \lam_w N \lam_{\min}r,
$$
where $\lam_{\min} = \min\{\lam_1,\ldots,\lam_m\}$, by the definition of
$\Wk(r)$.
Recall that $n$ can be chosen arbitrarily large, so $\lam_w$ can be
arbitrarily small, and we obtain that
$$
\ov{D}(\nu_a,\th') \ge c_2^{-1}\lam_{\min} N.
$$
Since $N\in \Nat$ is arbitrary, the proposition follows. \qed


\section{Proof of the recurrence lemma~\ref{lem-rec}}

Let $K\in
\left\{0,\dots ,m\right\}$ be the number of $i$ for which
  $\varphi _i\not\in \pi \mathbb{Q}$.
Without loss of generality we may assume the following:
   if $K\ge 1$ then $\varphi _1,\dots ,\varphi _K\not\in \pi\Q$.

\smallskip

We distinguish the following cases:

\begin{description}
  \item[A] $\varphi _i \in \pi\Q$ for all $i\le m$.
  \item[B] there exists $i$ such that $\varphi _i\not\in \pi \mathbb{Q}$
  and $\varepsilon _i=1$.
  \item[C] $K\geq 1$ and $\varepsilon _i=-1$ for all $i\le K$.
  \begin{description}
               \item[C1] there exist $i,j\leq K$ such that
               $\varphi _i-\varphi _j\not\in \pi \mathbb{Q}$.
               \item[C2] there exists $r_i\in \mathbb{Q}$ such that
               $\varphi _i=\varphi _1+r_i\pi $ for $1\leq i\leq K$.
                \begin{description}
                              \item[C2a] $K<m$ and there  exists $j\geq K+1$
                              such that $\varepsilon _j=-1$.
                              \item[C2b] $K<m$ and for all $j\geq K+1$
                              we have $\varepsilon _j=1$.
                              \item[C2c] $K=m$.
                            \end{description}
             \end{description}
\end{description}

Denote by $R_\phi$ the rotation through the angle $\phi$. We call it
an irrational rotation if $\phi\not\in \pi\Q$.
Consider the semigroup generated by $\Ok_i,\ i\le m$, which we denote by
$\Sk$. We begin with the following observation.

\smallskip

{\sc Claim.} {\em Either $\Sk$ is finite, or $\Sk$ contains an
irrational rotation.}

\smallskip

The semigroup $\Sk$ is clearly finite in Case A and contains an
irrational rotation in Case B. In Case C1 we have $\Ok_i\Ok_j =
R_{\phi_i-\phi_j}$, which is an irrational rotation. In Case C2a we
also have that $\Ok_i\Ok_j = R_{\phi_i-\phi_j}$ is an irrational
rotation, since $\phi\not\in \pi\Q$ and $\phi_j\in \pi\Q$. We claim
that in remaining Cases C2b and C2c the semigroup is finite. This
follows easily; then $\Sk$ is generated by one irrational reflection
and finitely many rational rotations.

\medskip

{\em Proof of Lemma~\ref{lem-rec} when $\Sk$ is finite.} A finite
semigroup of invertible transformations is necessarily a group. Let
$\Sk = \{s_1,\ldots,s_t\}$. By the definition of the semigroup $\Sk$
we have $s_i = \Ok_{w^{(i)}}$ for some $w^{(i)} \in \Ak^*$,
$i=1,\ldots,t$. For every $v\in \Ak^*$ we can find $\widehat{v}\in
\Ak^*$ such that $\Ok_{\widehat{v}} = \Ok_v^{-1}$. Fix $u =
j_1\ldots j_k$ from the statement of the lemma. Consider the
following finite word over the alphabet $\Ak$:
$$
\om := \tau_1\ldots \tau_t,\ \ \ \mbox{where}\ \ \tau_j = (w^{(j)}
u)\,  \widehat{(w^{(j)}u)},\ j=1,\ldots,t.
$$
Note that $\Ok_{\tau_j} = I$ (the identity).
By the definition of $\Om$, the sequence $\bi\in \Om$ contains $\om$
infinitely many times. Suppose that $\sig^\ell\bi \in [\om]$. Since
$\Ok_{\bi|\ell} \in \Sk$, there exists $w^{(j)}$ such that $\Ok_{w^{(j)}} =
\Ok_{\bi|\ell}^{-1}$. Then the occurrence of $u$ in $\tau_j$, the $j$th
factor of $\om$, will be at the position $n$ such that $\Ok_{\bi|n}= I$,
so we will have $\phi_{\bi|n} = 0 \in [0,\delta]$ and
$\eps_{\bi|n} = 1$, as desired.

\medskip

{\em Proof of Lemma~\ref{lem-rec} when $\Sk$ is infinite.}
By the claim above, there exists $w\in \Ak^*$ such that $\phi_w \not\in\pi
\Q$ and $\eps_w= 1$. Fix $u = j_1\ldots j_k$ from the statement of the lemma.
Let
$$
v := \left\{\begin{array}{ll} uu, & \mbox{if}\ \phi_u\not\in\pi\Q;\\
                            uuw, & \mbox{if}\ \phi_u\in\pi\Q.\end{array}
\right.
$$
Observe that $\phi_v\not\in\pi\Q$ and $\eps_v =1$. Let $v^k=v\ldots v$
(the word $v$ repeated $k$ times). 
 Since $\phi_v/\pi$ is irrational, there exists an $N$ such that
every orbit of $R_{\phi_v}$ of length $N$
contains a point in every subinterval
of $[0,2\pi)$ of length $\delta$. Put
$$
\om:= \left\{\begin{array}{ll} v^N, & \mbox{if}\ \eps_i =1,\ \forall\, i\le m;
\\ v^N j^* v^N, & \mbox{if}\ \exists\,j^*\ \mbox{such that}\ \eps_{j^*} = -1.
\end{array}
\right.
$$
By the definition of $\Om$, the sequence $\bi\in \Om$ contains $\om$
infinitely many times. Let $\ell\in \Nat$ be such
that $\sig^\ell\bi \in [\om]$.
Suppose first that $\eps_{\bi|\ell} = 1$.
Then we have, denoting the length of $v$ by $|v|$,
\begin{equation} \label{rot}
\sig^{\ell + k|v|}\bi \in [u],\ \ \ \ \phi_{\bi|(\ell+k|v|)} = \phi_{\bi|\ell} +
k\phi_v\ (\mbox{mod}\ 2\pi),\ \ \ \ \eps_{\bi|(\ell+k|v|)} =1,
\end{equation}
for $k=0,\ldots, N-1$.
By the choice of $N$, we
can find $k\in \{0,\ldots,N-1\}$ such that $\phi_{\bi|(\ell+k|v|)}\in
[0,\delta]$, then $n=\ell + k|v|$ will be as desired.
If $\eps_{\bi|\ell} = -1$, then we replace $\ell$ by $\ell^*:=  \ell +
N|v|+1$ in (\ref{rot}), that is, we consider the occurrences of $u$ in
the second factor $v^N$. The orientation will be switched by
$\Ok_{j^*}$ and we can find the desired $n$ analogously.
\qed


\section{Concluding remarks}

Consider the special case when the self-similar set $\Lam$ is of the
form
\begin{equation} \label{eq-ss2}
\Lam = \bigcup_{i=1}^m (\lam_i \Lam + b_i),\ \ \ b_i\in \R^2.
\end{equation}
In other words, the contracting similitudes have no rotations or reflections,
as for the four corner Cantor set.
Then the projection $\Lam^\th:=p_\th(\Lam)$ is itself a self-similar set on the
line:
$$
\Lam^\th = \bigcup_{i=1}^m (\lam_i \Lam^\th +
p_\th(b_i)),\ \ \mbox{ for}\  \th\in [0,\pi).
$$
Let $\Lam^\th_i = \lam_i \Lam^\th + p_\th(b_i)$. As above, $\nu$ is the natural
measure on $\Lam$. Let $\nu_\th$ be the natural measure on $\Lam^\th$,
so that $\nu_\th = \nu \circ p_\th^{-1}$.

\begin{corollary} \label{cor-pro}
Let $\Lam$ be a self-similar set of the form (\ref{eq-ss2}) that is not
on a line, such that $\sum_{i=1}^m \lam_i \le 1$.
If $\Lam$ satisfies the Open Set Condition condition, then
$$
\nu_\th(\Lam^\th_i \cap \Lam^\th_j) = 0,\ i\ne j,\ \ \ \mbox{for a.e.}
\ \th\in [0,\pi).
$$
\end{corollary}

{\em Proof.} Let $s>0$ be such that $\sum_{i=1}^m \lam_i^s = 1$.
By assumption, we have $s\le 1$. This number is known as the similarity
dimension of $\Lam$ (and also of $\Lam^\th$ for all $\th$). Suppose first
that $s=1$. Then we are in the situation covered by Theorem~\ref{th-main},
and $\nu$ is just the normalized restriction of $\Hk^1$ to $\Lam$.
Consider the product measure $\nu\times \Lk$, where $\Lk$ is the Lebesgue
measure on $[0, \pi)$. Theorem~\ref{th-main} implies that
$$
(\nu\times \Lk)\{(x,\th)\in \Lam\times [0,\pi):\ \exists\,y\in \Lam,\ y\ne x,
\ p_\th(x) = p_\th(y)\} = 0.
$$
By Fubini's Theorem, it follows that for $\Lk$ a.e.\ $\th$,
for $\nu_\th$ a.e.\ $z\in L^\th$, we have that $p_\th^{-1}(z)$ is a
single point. This proves the desired statement, in view of the fact that
$\nu(\Lam_i \cap \Lam_j) = 0$ for $\Lam$ satisfying the 
Open Set Condition.

In the case when $s<1$ we can use \cite[Proposition 1.3]{PSS}, which
implies that the packing measure $\Pk^s(\Lam^\th)$ is positive and finite
for $\Lk$ a.e.\ $\th$. By self-similarity and the properties of $\Pk^s$
(translation invariance and scaling), we
have $\Pk^s(\Lam^\th_i\cap \Lam^\th_j) = 0$ for $i\ne j$.
Then we use \cite[Corollary 2.2]{PSS},
which implies that $\nu_\th$ is the normalized restriction of
$\Pk^s$ to $\Lam^\th$, to complete the proof. \qed

\medskip

\noindent {\bf Remark.}
In \cite[Proposition 2]{BG} it is claimed that if a self-similar set
$\Kk = \bigcup_{i=1}^m \Kk_i$ in $\R^d$ has the
Hausdorff dimension equal to the similarity dimension, then the natural
measure of the ``overlap set'' $\bigcup_{i\ne j}(\Kk_i\cap \Kk_j)$
is zero. This would imply Corollary~\ref{cor-pro},
since the Hausdorff dimension of
$\Lam^\th$ equals $s$ for $\Lk$ a.e.\ $\th$ by Marstrand's Projection
Theorem. Unfortunately, the proof in \cite{BG} contains an error,
and it is still  unknown whether
the result holds [C. Bandt, personal communication].
(It should be noted that \cite[Proposition 2]{BG}
was not used anywhere in \cite{BG}.)

\medskip

\noindent
{\bf Acknowledgment.} We are grateful to M. Cs\"ornyei, E. J\"arvenp\"a\"a,
and M. J\"arvenp\"a\"a for helpful discussions. This work was done while
K. S. was visiting the University of Washington.


\end{document}